\documentclass[english,a4paper,12pt]{article}
\usepackage{babel}
\usepackage[ansinew]{inputenc}
\usepackage{graphicx}
\usepackage{verbatim}
\usepackage{amsmath}
\usepackage{amssymb}
\def\C{\mathbb C}
\def\Cd{\widehat{\mathbb C}}
\def\N{\mathbb N}

\newtheorem{theo}{Theorem}[section]
\newtheorem{ex}[theo]{Example}

\newtheorem{cor}[theo]{Corollary}

\title{A remark concerning normal families and shared values}
\author{Andreas Sauer*}
\begin{document}
\maketitle
\thispagestyle{empty}
\textit{\small * Hochschule Ruhr West, Duisburger Str.~100, 45479 M\"{u}lheim an der Ruhr,\\ Germany, andreas.sauer@hs-ruhrwest.de, ORCID: 0000-0003-0542-4745}\\
%\textit{\small ** Department of Mathematics Education, Kongju National University, Gongju 32588,\\South Korea, schweizer@kongju.ac.kr, ORCID: 0000-0002-3097-993X}}
\begin{abstract} We improve well-known results concerning normal families and shared values of meromorphic functions in the plane. In particular, we obtain two corollaries concerning meromorphic functions $f \colon {\mathbb C} \to {\widehat{\mathbb C}}$: i) If $f$ shares a non-zero finite value with $f'$, and such that $f'$ is bounded on the preimages of $f$ for a second value, then $f$ is normal. ii) If $f$ shares two finite values with $f'$, then $f$ and $f'$ are normal.
\end{abstract}\vskip.3cm
\text{\small 2020 Mathematics Subject Classification: Primary: 30D45, Secondary: 30D35\\}\\
\renewcommand{\thefootnote}{}
\footnotetext{\hspace*{-.51cm} %
}
%\noindent
\section{Introduction}
Let $D \subset \C$ be a domain and denote by $\cal M$ the set of all meromorphic functions on $D$ endowed with the topology of locally uniform convergence. A subset $\cal F \subset \cal M$ is called a {\it normal family} if $\cal F$ is pre-compact in $\cal M$. For the theory of normal families of meromorphic functions and a brief introduction to Nevanlinna theory we refer to \cite{Schi}. In particular, we will use Marty's criterion on several occasions (\cite{Schi}, section 3.3). For more details concerning Nevanlinna theory we refer to \cite{Ya}.\\

The criteria that we consider, use three typical assumptions for meromorphic functions $f \colon D \to \Cd$, that are related to each other, and can be described as follows: The strongest such assumption is that for a complex number $a \in \C$ it holds $f(z)=a\Leftrightarrow f'(z)=a$ for all $z \in D$, when $a$ is called a {\it shared} value of $f$ and $f'$. A weaker requirement is that for some $a \in \C$ it holds $f(z)=a \Rightarrow f'(z)=a$ for all $z \in D$, when in recent years this has often been referred to as a {\it partially} shared value of $f$ and $f'$. An even weaker assumption is that for $a \in \C$ it holds $f(z)=a \Rightarrow |f'(z)| \le K$ for all $z \in D$ with a constant $K$. We will call a value $a$ with this property a {\it value with bounded derivative} for $f$.\\

The amount of publications on the connection between normal families of meromorphic functions and these three assumptions for a certain number of values for each $f \in \cal F$ is vast, and we make no attempt to give a complete overview. But we want to mention some important references: An easy consequence of the results of Lappan in \cite{La} is, that five values with bounded derivative imply normality. Schwick proved in \cite{Schw} that the same conclusion holds for three shared values. Later Pang and Zalcman reduced this in \cite{PaZa} to two shared values. In \cite{XuFa} Xu and Fang proved that three partially shared values imply normality, and in \cite{ChFaZa} Chang, Fang and Zalcman gave a proof for one non-zero shared value and one further partially shared value.\\

We will continue these investigations by improving two of the above mentioned results.\\

First we show in Theorem \ref{more_general} that one non-zero shared value and one further value with bounded derivative imply normality. To achieve this, we give a refined version of the proofs in \cite{PaZa} and \cite{ChFaZa}.\\

Secondly, we complement the results in \cite{PaZa} by showing, that two shared values do not only imply normality of the family $\cal F$ itself, but also the normality of the family of all derivatives of functions in $\cal F$.\\

In most papers on this topic the lemmata named after Pang and Zalcman are employed (more or less including \cite{La}), and the current paper is no exception to this. We refer to \cite{PaZa2}.

\section{One non-zero shared value and one value with bounded derivative}
The main result of this section is the following theorem.
\begin{theo} \label{more_general}
Let $\cal F$ be a family of meromorphic functions on a domain $D \subset \C$, $a \neq 0$ and $b \neq 0$ be complex numbers and $K \ge 1$. If for every $f \in \cal F$ and all $z \in D$
$$
f(z) = a \Leftrightarrow f'(z) = b \quad \text{and} \quad f(z) = 0 \Rightarrow |f'(z)| \le K, 
$$ 
then $\cal F$ is a normal family.
\end{theo}
Let $f \colon \C \to \Cd$ be a meromorphic function. Then $f$ is {\it normal}, if for every sequence $z_n \to \infty$ the family formed by the functions $f(z_n + z)$ is normal. This is equivalent to the boundedness of the spherical derivative $f^{\#}(z) = |f'(z)|/(1+|f(z)|^2)$ for all $z \in \C$. The property of normality for $f \colon \C \to \Cd$ is rather restrictive (see \cite{Er}), and is often an important argument in proofs of uniqueness theorems concerning $f$ and $f'$.\\

In analogy to Theorem 3 in \cite{PaZa} we get:
\begin{theo} \label{f_normal}
Let $f \colon \C \to \Cd$ be a meromorphic function and $a$ and $b \neq 0$ be distinct complex numbers. If $f'$ is bounded on $f^{-1}(a)$ and if $f$ and $f'$ share $b$, then $f$ is normal.
\end{theo}
\noindent
Theorem \ref{f_normal} immediately gives the following corollary, which also follows from \cite{ChFaZa}:

\begin{cor} \label{partial_normal}
Let $f \colon \C \to \Cd$ be a meromorphic function and $a$ and $b \neq 0$ be distinct complex numbers. If $f(z)=a \Rightarrow f'(z)=a$ for all $z \in \C$, i.e.~if $a$ is partially shared by $f$ and $f'$, and if $f$ and $f'$ share $b$, then $f$ is normal.
\end{cor}
Since Theorem \ref{f_normal} follows from Theorem \ref{more_general} by considering $g(z) = f(z) - a$ and the related families formed by the functions $g(z_n + z)$ with arbitrary $z_n \to \infty$ we only have to prove Theorem \ref{more_general}.\\

We give the following example:
\begin{ex} \label{example_tan} \rm We define the family $\cal F$ on an arbitrary domain $D \subset \C$, consisting of the functions
$$
f_{n,c}(z) = 1 + \sqrt{1+n} \tan \left( \frac{z+c}{\sqrt{1+n}} \right)
$$
with $n \in \N$ and $c \in \C$. A simple calculation shows
$$
f'_{n,c}(z) = 1 + \tan \left( \frac{z+c}{\sqrt{1+n}} \right)^2.
$$
It is easy to check that all $f_{n,c}$ share the value $1$ with $f'_{n,c}$ and that $f_{n,c}(z)=0$ implies $f'_{n,c}(z) = 1 + 1/(1+n)$, so that for all $f \in \cal F$ we have $f(z)=0 \Rightarrow |f'(z)| \le 3/2$. Also there exists no value other than $1$ that is partially shared by $f_{n,c}$ and $f'_{n,c}$ for all $n$ and $c$. ($f_{n,c}$ and $f'_{n,c}$ partially share the value $n+2$, but this is not a global value for all $f \in \cal F$.) Hence Theorem \ref{more_general} shows that $\cal F$ is a normal family, while the criteria in \cite{PaZa} and \cite{ChFaZa} would not.\\

Example \ref{example_tan} is only an illustration. As in most cases where a normal family $\cal F$ consists of elementary functions, it is easily possible to show that the spherical derivatives of all $f \in \cal F$ are uniformly bounded, so that no special criteria are needed to show normality.
\end{ex}

\noindent
Inspection of the proof of Theorem \ref{more_general} shows the following for analytic functions.
\begin{theo} \label{more_general_analytic}
Let $\cal F$ be a family of holomorphic functions on a domain $D \subset \C$, $a \neq 0$ and $b \neq 0$ be complex numbers and $K \ge 1$. If for every $f \in \cal F$ and all $z \in D$
$$
f'(z) = b \Rightarrow f(z) = a \quad \text{and} \quad f(z) = 0 \Rightarrow |f'(z)| \le K, 
$$ 
then $\cal F$ is a normal family.
\end{theo}

\noindent
This again gives a corollary.
\begin{theo} \label{f_normal_entire}
Let $f \colon \C \to \C$ be an entire function and $a$ and $b \neq 0$ be distinct complex numbers. If $f'$ is bounded on $f^{-1}(a)$ and if $f'(z) = b$ implies $f(z)=b$ for all $z \in \C$, then $f$ is normal.
\end{theo}
\noindent
{\it Proof of Theorem \ref{more_general}.} Suppose $\cal F$ is not a normal family. Then, by Lemma 1 in \cite{PaZa}, we have $f_n \in \cal F$, $z_n \to z_0 \in D$ and $\rho_n > 0$ with $\rho_n \to 0$, such that
$$
g_n(\zeta) := \frac{f_n(z_n + \rho_n \zeta)}{\rho_n} \to g(\zeta)
$$
locally uniformly with respect to the spherical metric, where $g \colon \C \to \Cd$ is a non-constant normal meromorphic function with $g^{\#}(0)=K+1$.\\

First we show that $g'$ cannot be constant. Suppose to the contrary that $g' \equiv A$ with $A \neq 0$, so that $g(\zeta) = A \zeta + B$ with $B \in \C$. Then $g$ has a zero at $\zeta_0 = -B/A$, so that there exists a sequence $\zeta_n \to \zeta_0$ with
$$
g_n(\zeta_n) = \frac{f_n(z_n + \rho_n \zeta_n)}{\rho_n} = 0, 
$$
hence $f_n(z_n + \rho_n \zeta_n) = 0$. By assumption it follows $|g'_n(\zeta_n)| = |f_n'(z_n + \rho_n \zeta_n)| \le K$. From $g'_n(\zeta_n) \to g'(\zeta_0) = A$ we get $|A| \le K$. But this implies $g^{\#}(0) \le |g'(0)| = |A| \le K$, in contradiction to  $g^{\#}(0)=K+1$.\\

We prove $g(\zeta) = \infty$ if and only if $g'(\zeta)=b$. This implies that $g$ has no poles, i.e.~$g$ is entire, and that $g'$ omits the value $b$. Under slightly stronger assumptions this was proved in \cite{ChFaZa}, Proof of Theorem 2, (ii) and (iii). For convenience we repeat the arguments given there, which can be applied unaltered under the present assumptions.\\

Suppose $g'(\zeta_0) = b$. Since $g'$ is not constant there exists a sequence $\zeta_n \to \zeta_0$ such that
$$
g'_n(\zeta_n) = f_n'(z_n + \rho_n \zeta_n) = b,
$$
and therefore by assumption
$$
g_n(\zeta_n) = \frac{f_n(z_n + \rho_n \zeta_n)}{\rho_n} = \frac{a}{\rho_n},
$$
so that $g(\zeta_0) = \lim_{n \to \infty} g_n(\zeta_n) = \infty$.\\

Now suppose $g(\zeta_0) = \infty$ of order $m$. Then $1/g(\zeta_0) = 0$ of order $m$, and since $1/g_n(\zeta) - \rho_n/a \to 1/g(\zeta)$ for large $n$ there exist $m$ solutions $\zeta_{n,1}, \ldots, \zeta_{n,m}$ (when counted according to multiplicities) of the equation $1/g_n(\zeta) - \rho_n/a = 0$ that converge to $\zeta_0$ with $n \to \infty$. For these $\zeta_{n,j}$ it follows
$$
g_n(\zeta_{n,j}) - \frac{a}{\rho_n} = \frac{f_n(z_n + \rho_n \zeta_{n,j})-a}{\rho_n}= 0.
$$
Hence $f_n(z_n + \rho_n \zeta_{n,j}) = a$, so that by assumption $f_n'(z_n + \rho_n \zeta_{n,j}) = b$. It follows $(1/g_n)'(\zeta_{n,j}) = -g_n'(\zeta_{n,j})/g_n^2(\zeta_{n,j}) = -b \rho_n^2/a^2 \neq 0$, so that all the zeros $\zeta_{n,j}$ of $1/g_n(\zeta) - \rho_n/a$ are simple, and therefore $m$ distinct points. Since $1/g$ is holomorphic at $\zeta_0$ and $1/g_n \to 1/g$ it follows that $(1/g_n)' = -g_n'/g_n^2 \to (1/g)'$ locally uniformly. We deduce that the sequence $- g_n'/g_n^2 + b \rho_n^2/a^2$ converges to $(1/g)'$ locally at $\zeta_0$ with at least $m$ zeros $\zeta_{n,j} \to \zeta_0$, so that according to Rouch{\'e}'s theorem $(1/g)'$ has a zero at $\zeta_0$ at least of order $m$, which is a contradiction.\\
%$$
%g'(\zeta_0) = \lim_{n \to \infty}  g'_n(\zeta_n) = \lim_{n \to \infty} f_n'(z_n + \rho_n \zeta_n) = b.
%$$

As mentioned above, it follows that $g$ is a normal entire function, such that $g'$ omits the value $b$. Since $g$ is normal, the order of $g$ is at most one, as follows from a classical result in \cite{ClHa}. Hence the order of $g'$ is at most one (see e.g.~Theorem 4.2 in \cite{Ya}). It follows that $g'$ has to be of the form
$$
g'(\zeta) = C_1 \text{e}^{\lambda \zeta} + b,
$$
with non-zero constants $C_1$ and $\lambda$, so that
$$
g(z) = \frac{C_1}{\lambda} \text{e}^{\lambda \zeta} + b\zeta + C_2.
$$
But then $g$ is not normal: $g$ has infinitely many zeros $\zeta_n \to \infty$, since otherwise
$$
\frac{-C_1 \text{e}^{\lambda \zeta}}{\lambda (b\zeta + C_2)} 
$$  
is a non-constant meromorphic function with one pole, no zeros and only finitely many $1$-points, contradicting Picard's theorem. We immediately get 
$$
g'(\zeta_n) = -\lambda(b\zeta_n + C_2) + b \to \infty,
$$
and hence $g^{\#}(\zeta_n) = |g'(\zeta_n)| \to \infty$, in contradiction to the normality of $g$. $\square$
\section{Two shared values and derivatives}
For a family of meromorphic functions $\cal F$ we will denote by $\cal F'$ the family of all derivatives of functions in $\cal F$. The only function that does not have a well-defined derivative is $f \equiv \infty$. We will assume that $f \equiv \infty$ is not in $\cal F$. Note that a single function cannot change the property of normality for a family $\cal F$.\\

The emphasis of the following theorem lies on the statement concerning $\cal F'$. It is well-known that for a normal family $\cal F$ of meromorphic functions, $\cal F'$ need not be normal (\cite{Schi}, Example 3.1.8).
\begin{theo} \label{two_shared_values_derivative}
Let $\cal F$ be a family of meromorphic functions on a domain $D \subset \C$, $a$, $b$, $c$ and $d$ be complex numbers with $a \neq b$. If for every $f \in \cal F$ and all $z \in D$
$$
f(z) = a \Leftrightarrow f'(z) = c \quad \text{and} \quad f(z) = b \Leftrightarrow f'(z) = d, 
$$ 
then $\cal F$ and $\cal F'$ are normal families.
\end{theo}

This gives the following corollary, which except for the statement concerning $f'$ is Theorem 3 in \cite{PaZa}.

\begin{cor} \label{derivative_normal}
Let $f \colon \C \to \Cd$ be a meromorphic function such that $f$ and $f'$ share two finite values, then $f$ and $f'$ are normal functions.
\end{cor}

We give an example that shows that the derivative of a normal meromorphic function need not be normal.

\begin{ex} \rm Consider
$$
f(z) := \frac{1}{\text{e}^z + \frac{1}{z^2}},
$$
so that
$$
f'(z) = \frac{ \frac{2}{z^3} - \text{e}^z}{\left( \text{e}^z + \frac{1}{z^2} \right)^2}.
$$
It is easy to see that $f$ is normal: $\text{e}^z$ is normal, so that $\text{e}^z + 1/z^2$ is normal.\\

Let $z_n \to \infty$ be a sequence of zeros of $f'$, which means $\text{e}^{z_n} = 2/z_n^3$. Direct calculation shows
$$
f''(z_n) = - \frac{2(z_n + 3) z_n^2}{(z_n + 2)^2} \to \infty,
$$
hence $(f')^{\#}(z_n) = |f''(z_n)| \to \infty$, and $f'$ is not normal.
\end{ex}

{\it Proof of Theorem \ref{two_shared_values_derivative}.} Theorem 2 in \cite{PaZa} shows that $\cal F$ is normal. (See \cite{ChFaZa}, Theorem 2 for a shorter proof, or simply apply Theorem \ref{more_general} of the current paper.)\\ 

Suppose $\cal F'$ is not normal. Then there exists a sequence $f'_n \in \cal F'$ and a sequence $z_n \to z_0 \in D$ with $(f'_n)^{\#}(z_n) \to \infty$, so that $f_n'(z)$ is not normal at $z_0$. Let $f_n$ be a corresponding sequence of primitive functions contained in $\cal F$. Since $\cal F$ is a normal family we can assume by passing to a subsequence, that $f_n$ converges compactly in $D$ to a meromorphic function $F \colon D \to \Cd$, where $F \equiv \infty$ is possible.\\

According to Zalcman's Lemma, there exist sequences $w_n \to z_0$ and $\rho_n \to 0$ such that a further subsequence
$$
\psi_n(\zeta) := f_n'(w_n + \rho_n \zeta)
$$
converges locally uniformly in the $\zeta$-plane to $\Psi(\zeta)$, which is a non-constant normal meromorphic function $\Psi \colon \C \to \Cd$. Since $w_n + \rho_n \zeta \to z_0$ locally uniformly in the $\zeta$-plane, it is immediate that $f_n(w_n + \rho_n \zeta) \to F(z_0)$ locally uniformly.\\

We show $F(z_0) = a$ or $F(z_0) = b$.\\

Suppose $\Psi$ omits the values $c$ and $d$. Note that the assumptions imply $c \neq d$. Because all poles of $\Psi$ are multiple, the second fundamental theorem of Nevanlinna theory (Theorems 1.5 and 1.6 in \cite{Ya}) shows 
\begin{align*}
T(r,\Psi) & \le \overline{N}(r,\Psi) + \overline{N}(r,c,\Psi) + \overline{N}(r,d,\Psi) + O(\log r) \\
& \le \frac{1}{2} N(r,\Psi) + O(\log r)\\ 
& \le \frac{1}{2} T(r,\Psi) + O(\log r),
\end{align*}
so that $T(r,\Psi) = O(\log r)$. Note that since $\Psi$ is normal, its order is at most $2$, so the error term $O(\log r)$ is justified. Hence $\Psi$ is a non-constant rational function. But then $\Psi$ can omit at most one value on $\C$, namely its value at $\zeta = \infty$, a contradiction.\\

Thus we can assume, without loss of generality, that there exists $\zeta_0 \in \C$ such that $\Psi(\zeta_0) = c$. Since $\Psi$ is non-constant there exists a sequence $\zeta_n \to \zeta_0$, such that $\psi_n(\zeta_n) = c$. By the assumptions it follows that $f_n(w_n + \rho_n \zeta_n) = a$, which implies $F(z_0)=a$.\\

Let $K$ be a compact disk in $D$ that contains $z_0$. Since $f_n$ is normal, there exists $M > 0$ such that $f_n^{\#}(z) \le M$ for all $z \in K$. $\Psi$ is non-constant, so there exists $\zeta^*$ such that $|\Psi(\zeta^*)| > 2 M (1 + |a|^2)$. Hence there exists a sequence $\zeta_n \to \zeta^*$ such that $|\psi_n(\zeta_n)| \ge 2 M (1 + |a|^2)$. Since $f_n(w_n + \rho_n \zeta) \to a$ locally uniformly in the whole $\zeta$-plane we conclude $|f_n(w_n + \rho_n \zeta_n)| \le |a| + \varepsilon$ for arbitrary $\varepsilon > 0$ and $n$ large enough. But then
$$
f_n^{\#}(w_n + \rho_n \zeta_n) \ge \frac{2 M (1 + |a|^2)}{1 + (|a| + \varepsilon)^2} > M
$$
for $\varepsilon$ small enough. Since for large $n$ we have $w_n + \rho_n \zeta_n \in K$ we obtain a contradiction. Hence $\cal F'$ is normal. $\square$\\
\vskip.3cm
\noindent
\noindent
{\bf Compliance with Ethical Standards:}\\
\noindent
Conflict of Interest: The Author declares that he has no conflict of interest.\\
\noindent
Ethical approval: This article does not contain any studies with human participants or animals performed by any of the authors.


\begin{thebibliography}{9999999}
%
\bibitem[ChFaZa]{ChFaZa} J.~Chang, M.~Fang and L.~Zalcman, Normality and fixed-points of meromorphic functions. {\it Ark.~Mat.} {\bf 43}, No. 2, (2005) 307--321. 
%
\bibitem[ClHa]{ClHa} J.~Clunie and W.K.~Hayman, The spherical derivative of integral and meromorphic functions. {\it Comment. Math. Helv.}
{\bf 40} (1966) 117--148.
%
\bibitem[Er]{Er} A.~Eremenko, Normal holomorphic curves from parabolic regions to projective spaces. {arXiv:0710.1281}, (2007).
%
\bibitem[La]{La} P.~Lappan, A criterion for a meromorphic function to be normal. {\it Comment.~Math.~Helv.} {\bf 49}, (1974) 492--495. 
%
\bibitem[PaZa]{PaZa} X.~Pang and L.~Zalcman, Normality and shared values. {\it Ark.~Mat.} {\bf 38}, No. 1, (2000) 171--182.
%
\bibitem[PaZa2]{PaZa2} X.~Pang and L.~Zalcman, Normal families and shared values. {\it Bull.~Lond.~Math.~Soc.} {\bf 32}, No. 3, (2000) 325--331. 
%
\bibitem[Schi]{Schi} J.L.~Schiff, Normal families. {\it Springer-Verlag}, 1993.
%
\bibitem[Schw]{Schw} W.~Schwick, Sharing values and normality. {\it Arch.~Math.} {\bf 59}, No. 1, (1992) 50--54. 
%
\bibitem[XuFa]{XuFa} Y.~Xu and M.~Fang, A note on some results of Schwick. {\it Bull.~Malays.~Math.~Sci.~Soc. (2)}, {\bf 27}, No.~1, (2004) 1--8.
%
\bibitem[Ya]{Ya} L.~Yang, Value distribution theory. {\it Berlin: Springer-Verlag. Beijing: Science Press}, (1993). 
%
\end{thebibliography}
\end{document}